\documentclass[a4paper,11pt]{article} 
\usepackage{amsmath, amssymb, amsthm, shadow,epsfig, color}
\usepackage{thmtools}

\newcommand{\form}{\mathcal E}
\newcommand{\dom}{\mathcal F}

\newcommand{\dis}{\displaystyle}
\newcommand{\rd}{{\mathbb R}^d}
\newcommand{\vareps}{\varepsilon}
\newcommand{\wt}{\widetilde}

\numberwithin{equation}{section}

\declaretheoremstyle[headfont=\normalfont]{normalhead}

\declaretheorem[numberwithin=section]{theorem}

\declaretheorem[numberlike=theorem]{lemma}
\declaretheorem[numberlike=theorem]{proposition}

\declaretheorem[numberlike=theorem, style=normalhead]{definition}
\declaretheorem[style=normalhead, numberlike=theorem]{remark}

\makeatletter
\newcommand{\subjclass}[2][2010]{%
  \let\@oldtitle\@title%
  \gdef\@title{\@oldtitle\footnotetext{#1 \emph{Mathematics subject classification.} #2}}%
}
\newcommand{\keywords}[1]{%
  \let\@@oldtitle\@title%
  \gdef\@title{\@@oldtitle\footnotetext{\emph{Key words and phrases.} #1.}}%
}
\makeatother

\allowdisplaybreaks[4]

\title{\sf On  instability of global path properties of symmetric Dirichlet forms  
under Mosco-convergence}

\author{Kohei Suzuki \and Toshihiro Uemura}
        

\begin{document}
\date{}
\maketitle
\begin{abstract}
We give sufficient conditions for Mosco convergences for the following three cases: symmetric locally uniformly elliptic diffusions, 
symmetric L\'evy processes, and symmetric jump processes in terms of 
the $L^1(\rd;dx)$-local convergence of the (elliptic) coefficients, the characteristic exponents and the jump density functions,
respectively. 
We stress that the global path properties of the corresponding Markov processes 
such as recurrence/transience, and conservativeness/explosion 
are not preserved under Mosco convergences and we give several examples where such situations indeed happen.

\end{abstract}

\section{Introduction}

In the present paper, we are concerned with Mosco convergences of the following 
three types of the Dirichlet forms: symmetric strong local Dirichlet forms satisfying the locally uniformly 
elliptic conditions, symmetric translation invariant Dirichlet forms, and symmetric jump-type Dirichlet forms. 
We give sufficient conditions for the Mosco convergences in the above three cases in terms of the $L^1$-local convergence of the corresponding coefficients, and show instability of global 
path properties under the Mosco convergences such as recurrence or transience, and conservativeness or explosion by giving several examples.

We find that the Mosco convergences follow from quite mild assumptions  (see {\bf Assumption (A)}, {\bf (B)} and {\bf (C)}), which are essentially $L^1(\rd;dx)$-local convergences of the corresponding coefficients, 
which are diffusion coefficients,  L\'evy exponents and jump densities. 
Here $dx$ denotes the Lebesgue measure on $\rd$. Hereafter we fix our state space to $(\rd, dx)$ and we write $L^p(\rd)$ (or $L^p$) shortly for $L^p(\rd;dx)$ $(1 \le p \le \infty)$.
Since the $L^1$-local convergence is one of the weakest notions of 
convergences, our results mean that the weakest convergence of the coefficients implies 
the Mosco convergences.

The Mosco convergence is a notion of convergences of closed forms on Hilbert spaces (see Definition 
\ref{defi: Mosco}), which was introduced by U. Mosco \cite{M67}, originally to study the 
approximations of some variational inequalities. In \cite{M94}, he showed that the Mosco convergence 
is equivalent to the strong convergences of the corresponding semigroups and resolvents. The 
strong  convergence of semigroups implies the convergence of finite-dimensional distributions of the 
corresponding Markov processes when closed forms in question are regular Dirichlet forms. 
For such reasons, the Mosco convergence has been used to show the weak convergence of stochastic 
processes in the theory of Markov processes (see {\it e.g.} \cite{U95, KU97, S98, Kl05, Km06} and 
references therein). In \cite{KS03}, Kuwae and Shioya generalized the notion of the Mosco convergence,
now is called the {\it Mosco-Kuwae-Shioya convergence}, as the basic $L^2$-spaces can 
change, while Hino considered the non-symmetric version of the Mosco convergence in \cite{H98}.
Although both generalizations are quite important, in the present paper, we consider only symmetric cases
and we fix a basic $L^2$-space as $L^2(\rd)$.

Our another aim is to show that the Mosco convergence of Dirichlet forms does not preserve any 
global path properties for the corresponding processes of the Dirichlet forms in any respect.
As stated above, the Mosco convergence is equivalent to the strong convergence of the corresponding semigroups, which implies only the convergences of finite-dimensional distributions
of the corresponding Markov processes. Thus it is easy to imagine that global properties such as recurrence/transience and conservativeness/explosion are not preserved under 
the Mosco convergence. It seems, however, that those studies how to construct such examples concretely have not been investigated.

In this paper, we construct several examples whose global properties such as recurrence/transience and conservativeness/explosion are not preserved under the Mosco convergence.
In constructing such examples, we use the results about sufficient conditions for Mosco convergences as explained in the second paragraph in this introduction.


To be more precise, let us first consider symmetric strongly local Dirichlet forms having  
the locally uniformly elliptic coefficients. Namely let $A_n(x)=(a^n_{ij}(x))$
be a sequence of $d\times d$ symmetric matrix valued functions  satisfying the following conditions: 

\medskip
\noindent
\underline{\bf Assumption (A):}
\begin{itemize}
\item[\sf (A1)] For any compact set $K\subset \rd$, there exists $\lambda=\lambda(K)>0$ so that 
for all $n\in{\mathbb N}$, 
$$
 \lambda |\xi|^2 \le \langle A_n(x) \xi, \xi \rangle  \le \lambda^{-1} |\xi|^2, \quad 
\forall \, x\in K, \ \forall \, \xi \in \rd.
$$
\item[\sf (A2)] For any compact set $K$, 
	\begin{align*}
	\int_K \|A_n(x)-A(x)\|dx \to 0 \quad (n \to \infty),
	\end{align*}
where $\dis \|A_n(x)-A(x)\|^2 := \sum_{i=1}^d\sum_{j=1}^d\bigl(a^n_{ij}(x)-a_{ij}(x)\bigr)^2, \ x\in \rd$.
\end{itemize}

Then, consider a sequence of symmetric quadratic forms
$$
\form^n(u,v)=\int_{\rd} \langle A_n(x)\nabla u(x), \nabla v(x)\rangle dx
$$
and
$$
\form(u,v)=\int_{\rd} \langle A(x)\nabla u(x), \nabla v(x)\rangle dx
$$
for $u$ and $v$ in $C_0^{\infty}(\rd)$, where $C_0^{\infty}(\rd)$ is the set of infinitely differentiable functions defined 
on $\rd$ with compact support. 
Under {\bf Assumption (A)}, it is known that 
$(\form^n, C_0^{\infty}(\rd))$ and $(\form, C_0^{\infty}(\rd))$ are Markovian closable forms on 
$L^2(\rd)$. They become regular symmetric Dirichlet forms $(\form^n, \dom^n)$ and $(\form, \dom)$
on $L^2(\rd)$.  Our first result is the following:

\begin{theorem} \label{thm: diffusion}
Suppose that Assumption {\sf (A)} holds.  Then the Dirichlet forms $(\form^n, \dom^n)$ converges 
 to $(\form,\dom)$ on $L^2(\rd)$ in the sense of Mosco.
\end{theorem} 
\begin{remark} \normalfont
\begin{itemize}
\item[{\sf (1)}] Mosco \cite{M94} gave several examples for which the Dirichlet forms converges 
in his sense. According to our diffusion forms example, assuming the convergence of the elliptic 
coefficients locally in $L^1(\rd)$, he have shown the $\Gamma$-convergence,  which is a bit weaker than his 
convergence.  He claimed that the convergence is indeed his convergence when the so-called 
``compactly imbedded" is held, which is a bit harder to show in general.

\item[{\sf (2)}]  In \cite{H98}, Hino has given several equivalent conditions in order that the semigroups converge 
strongly in $L^2$ corresponding to semi-Dirichlet forms not necessarily symmetric ones and also gave 
sufficient conditions for which the convergence holds.  In the case of diffusion type forms, his conditions 
are similar to that of ours (see \cite[Example 4.3]{H98}).

\item[{\sf (3)}] In \cite{LZ93} and \cite{RZ97}, they studied the convergence of quadratic forms under the
uniformly elliptic condition and obtained the weak convergence
of corresponding Markov processes.
In Theorem \ref{thm: diffusion} above, we only assume {\it the locally} uniformly elliptic condition. 
\end{itemize}
\end{remark}

We now consider the convergence of symmetric L\'evy processes. Let $\{\varphi_n\}$ be a sequence of the characteristic functions defined by symmetric 
convolution semigroups $\{\nu_t^n, t>0\}_{n\in{\mathbb N}}$: 
$$
e^{-t\varphi_n(x)}:= \hat{\nu}_t^n(x) \, \Bigl(=\int_{\rd} e^{i\langle x, y\rangle} \nu_t^n(dy) \, 
\Bigr), \quad x\in \rd.
$$
According to the L\'evy-Khinchin formula(\cite{FOT11}), we have the following characterization of 
$\nu_t^n$:
\begin{equation}
\dis \varphi_n(x) =\frac 12 \langle S_nx,x\rangle +\int_{\rd} \bigl(1-\cos(\langle x,y\rangle) n_n(dy), 
\label{LK-formula}
\end{equation}
\begin{equation}
\hspace*{-4cm} \mbox{where} \ S_n \mbox{ is a non-negative definite symmetric } d\times d 
\mbox{ matrix \ and} 
\label{LK-1}
\end{equation}
\begin{equation}
n_n(dy) \mbox{ is a symmetric Borel measure on }  \rd\setminus\{0\} \mbox{ so that}
\dis{
\int_{\rd\setminus\{0\}} \frac{|x|^2}{1+|x|^2} n_n(dx)<\infty.
}
\label{LK-2}
\end{equation}

\bigskip
We consider the following condition:

\bigskip
\noindent
\underline{\bf Assumption (B)}: \ $\varphi_n$ converges to a function $\varphi$ locally in $L^{1}(\rd)$.

\bigskip
Under the assumption, we find that $\varphi$ is also the characteristic function of a symmetric 
convolution semigroup $\{\nu_t, t>0\}$. 
Moreover the corresponding quadratic forms 
\begin{align*}
\form^n(u,v) & = \int_{\rd} \widehat{u}(x) \overline{\widehat{v}}(x) \varphi_n(x)dx \\
\dom^n &= \Bigl\{ u \in L^2(\rd) : \ \int_{\rd} \bigl| \widehat{u}(x)\bigr|^2 \varphi_n(x)dx <\infty 
\Bigr\}
\end{align*}
and 
\begin{align*}
\form(u,v) & = \int_{\rd} \widehat{u}(x) \overline{\widehat{v}}(x) \varphi(x)dx \\
\dom &= \Bigl\{ u \in L^2(\rd) : \ \int_{\rd} \bigl| \widehat{u}(x)\bigr|^2 \varphi(x)dx <\infty 
\Bigr\}
\end{align*}
are symmetric translation invariant Dirichlet forms on $L^2(\rd)$. We  show that 
$(\form^n, \dom^n)$ converges to $(\form, \dom)$ in the sense of Mosco under {\bf Assumption (B)}:
\begin{theorem} Assume that Assumption {\bf (B)} holds. Then  
$(\form^n, \dom^n)$ converges to $(\form, \dom)$ in the  sense of Mosco. 
\label{Mosco:Levy}
\end{theorem}
The point is that we only assume the convergence {\it locally} in 
$L^1(\rd)$ of the respective quantities and no other further assumptions are needed.

\smallskip
We next consider the convergence of symmetric jump-type Dirichlet forms. Let $\wt{J}(x,y)$ be a non-negative symmetric 
function on $\rd\times \rd \setminus {\sf diag}$ satisfying 
\begin{align} \label{eq: Loc1}
x\mapsto \int_{y\not=x} \Bigl(1\wedge d(x,y)^2\Bigr) \wt{J}(x,y)dy \in L^1_{\sf loc}(\rd).
\end{align}
Here {\sf diag} means that the diagonal set: ${\sf diag}=\{(x,x) : \ x\in \rd \}$. 
Consider the following quadratic form $\wt\form$ on $L^2(\rd)$: 
$$
\wt{\form}(u,v) = \frac{1}{2}\iint_{x\not=y} \bigl(u(x)-u(y)\bigr)
\bigl(v(x)-v(y)\bigr)\wt{J}(x,y)dxdy
$$
for functions $u,v\in C_0^{\sf lip}(\rd)$. Here $C_0^{\sf lip}(\rd)$ is  the set of all Lipschitz continuous functions on $\rd$ with 
compact support.
Under the condition \eqref{eq: Loc1}, it is also known that 
$(\form, C_0^{\sf lip}(\rd))$  is a closable Markovian symmetric form on $L^2(\rd)$. 
Then the smallest closed extension $(\form,\dom)$ is a regular Dirichlet form.  

Now take $J_n(x,y)$ and $J(x,y)$  non-negative symmetric functions on $\rd\times \rd \setminus {\sf diag}$ 
satisfying \eqref{eq: Loc1} in place of $\wt J(x,y)$ and then consider regular symmetric jump-type Dirichlet forms as follows:
$$
\left\{
\begin{array}{rcl}
\form^n(u,v) &=& \dis{\frac{1}{2}\iint_{x\not=y} \bigl(u(x)-u(y)\bigr)
\bigl(v(x)-v(y)\bigr)J_n(x,y)dxdy, } \\
& & \vspace{-10pt} \\
\dom^n &=&  \overline{C_0^{\sf lip}(\rd)}^{\sqrt{\form_1^n}},
\end{array}
\right.$$
and 
$$
\left\{
\begin{array}{rcl}
\form(u,v) &=& \dis \frac{1}{2}\iint_{x\not=y} \bigl(u(x)-u(y)\bigr)
\bigl(v(x)-v(y)\bigr)J(x,y)dxdy, \\
& & \vspace{-10pt} \\
\dom &=&  \overline{C_0^{\sf lip}(\rd)}^{\sqrt{\form_1}},
\end{array}
\right.
$$
where $\mathcal E_1(u,v)=\mathcal E(u,v)+(u,v)_{L^2(\rd)}$.
\medskip
We make the following assumption.

\bigskip
\noindent
\underline{{\bf Assumption (C):}} 
\begin{itemize}
\item[{\sf (i)}] $J_n(x,y) \le \wt J(x,y)$ for $dx\otimes dy\mbox{-a.e.} \ (x,y) \in \rd\times \rd \setminus 
{\sf diag}$ and $\forall n\in{\mathbb N}$.

\item[{\sf (ii)}]  $\{J_n(x,y)\}$ 
converges to $J(x,y)$ {\it locally} in $L^1(\rd{\times}\rd\setminus {\sf diag};dx \otimes dy)$.
\end{itemize}

\smallskip
\begin{theorem} \label{thm: JtoJ}  Assume {\bf (C)}. Then $(\form^n, \dom^n)$ converges to $(\form, \dom)$ in the sense of Mosco. 
\end{theorem}

From now on, by using the above theorems, we construct several examples whose global path properties are not preserved under 
the Mosco convergence. We first consider the instability of conservativeness/explosion of the symmetric 
diffusion processes.  Under the same settings of Theorem \ref{thm: diffusion}, let the diffusion coefficients 
be diagonal $A_n(x)=a_n(x)I$, where $I$ denotes the identity of $d \times d$ matrices. Then we have
the following result:
\begin{proposition} \label{prop: Cons/Exp}
The following results hold:
\begin{itemize}
\item[\sf (i)] {\rm (explosive ones to conservative one)}
If we set 
$$
\alpha_n(x)= (2+|x|)^{2}\Bigl(\log(2+|x|)\Bigr)^{1+1/n}, \quad  
\alpha(x) =  (2+|x|)^{2}\Bigl(\log(2+|x|)\Bigr)
$$ 
for $n\in{\mathbb N}$,  then $(\form^n,\dom^n)$ is explosive for any $n$ and converges in  the sense of Mosco to
the conservative Dirichlet form $(\form, \dom)$.

\item[\sf (ii)] {\rm (conservative ones to explosive one)}
If we set 
$$
\alpha_n(x)= (2+|x|)^{2-1/n}\Bigl(\log(2+|x|)\Bigr)^{2}, \quad  
\alpha(x) =  (2+|x|)^{2}\Bigl(\log(2+|x|)\Bigr)^{2}
$$ 
for $n \in \mathbb N$,  then $(\form^n,\dom^n)$ is conservative for any $n$ and converges in the sense of Mosco to the explosive Dirichlet 
form $(\form, \dom)$.
\end{itemize}
\end{proposition}

We now consider the instability of recurrence/transience  of the symmetric L\'evy processes. Let $\alpha$ 
and $\alpha_n$ be measurable functions on $[0,\infty)$ satisfying that there exist positive constants 
$\underline{\alpha}$ and $\overline{\alpha}$ so that 
$$
0< \underline{\alpha} \le \alpha_n(t) \le \overline{\alpha} <2, \quad \text{a.e. } t\in [0,\infty)
$$
and define L\'evy measures on $\rd$ as follows: 
$$
n_n(dx)=|x|^{-d-\alpha_n(|x|)}dx, \quad n(dx)=|x|^{-d-\alpha(|x|)}dx.
$$
Then the corresponding characteristic (L\'evy) exponents are given by 
$$
\varphi_n(x)=\int_{\rd} \Bigl(1-\cos (x \xi)\Bigr) n_n(d\xi), \quad 
\varphi(x)=\int_{\rd} \Bigl(1-\cos (x \xi)\Bigr) n(d\xi),
$$
respectively. Then the following Proposition holds: 

\begin{proposition} \label{prop: Rec/Tra} Let $n_n$ and $n$ be as above. Assume $d=1$.  Then the following results hold.
\begin{itemize}
\item[\sf (i)] {\rm (recurrent ones to transient one)}
If we set 
$$
\alpha_n(u)=1+1/n-\Bigl(\log(u+ e^2)\Bigr)^{-1/2}, \quad  \alpha(u) = 1-\Bigl(\log(u+e^2)\Bigr)^{-1/2}
$$ 
for $u\ge 0$ and $n \in \mathbb N$, then $(\form^n,\dom^n)$ is recurrent for any $n$ and converges in the sense of Mosco to the transient 
Dirichlet form $(\form, \dom)$.

\item[\sf (ii)] {\rm (transient ones to recurrent one)}
If we set 
$$
\alpha_n(u)= 1- \Bigl(\log \bigl(u+e^2\bigr)\Bigr)^{-(1-1/n)}, \quad  
\alpha(u) = 1- \Bigl(\log \bigl(u+e^2\bigr)\Bigr)^{-1}
$$ 
for $u\ge 0$ and $n \in \mathbb N$, then $(\form^n,\dom^n)$ is transient for any $n$ and converges in the sense of Mosco to the recurrent 
Dirichlet form $(\form, \dom)$. 
\end{itemize}
\end{proposition}
The point is the sharp criterion of the recurrence/transience for the 
stable-type processes (see {\it e.g.} Theorem 3.3 in \cite{U04} and Theorem \ref{thm: Rec} in Appendix in the 
present paper).

\begin{remark} \rm  We can give a rather simple example for which the symmetric Dirichlet forms corresponding 
to transient symmetric L\'evy processes converge to the symmetric Dirichlet form corresponding to a 
recurrent symmetric L\'evy process in the sense of Mosco: 

Assume $d=2$.  Consider a function $\varphi_n(x):=|x|^{2-1/n}, \ x\in{\mathbb R}^2$ for each $n$. 
Then $\varphi_n$ defines the characteristic exponent associated with a transient symmetric $(2-1/n)$-stable 
process on ${\mathbb R}^2$. Clearly $\varphi_n(x)$ converges to $\varphi(x):=|x|^2$ for all $x$ and 
the limit $\varphi(x)$ is the characteristic exponent associated with a $2$-dimensional Brownian motion that 
is recurrent. Note that this example shows not only the instability of (global) path properties but also the 
instability of path types. Namely, the jump processes converge to the diffusion process. 
We will discuss such instability of path types in a forthcoming paper.
\end{remark}

\medskip
We finally consider the instability of recurrence/transience of symmetric jump processes. 
Let $\alpha$ and $\alpha_n$ be measurable functions on $[0,\infty)$ satisfying that there exist positive 
constants $\underline{\alpha}$ and $\overline{\alpha}$ so that 
$$
0< \underline{\alpha} \le \alpha_n(u) \le \overline{\alpha} <2, \quad {\rm a. e. } \ 
u\in [0,\infty).
$$
Let $c(x)$ be a measurable function on $\rd$ satisfying that there exist $0<c<C<\infty$ so that 
$c\le c(x)\le C$ for all $x\in \rd$. We consider the following jump kernels:
$$
J_n(x,y)=(c(x)+1)|x-y|^{-d-\alpha_n(|x-y|)}, \quad x,y \in \rd, \ x\not=y.
$$
and
$$
J(x,y)=(c(x)+1)|x-y|^{-d-\alpha(|x-y|)}, \quad x,y \in \rd, \ x\not=y.
$$
We note that the corresponding jump processes are not necessarily L\'evy processes because $c(x)$ is not necessarily translation-invariant.
Even in this case,  we have the following result similar to Proposition \ref{prop: Rec/Tra}:
\begin{proposition} \label{prop: Rec/Tra_J}Let $J_n$ and $J$ be as above. Assume $d=1$.  Then the following results hold.
\begin{itemize}
\item[\sf (i)] {\rm (recurrent ones to transient one)}
If we set 
$$
\alpha_n(u)=1+1/n-\Bigl(\log(u+ e^2)\Bigr)^{-1/2}, \quad  \alpha(u) = 1-\Bigl(\log(u+e^2)\Bigr)^{-1/2}
$$ for $u\ge 0$ and $n \in \mathbb N$, then $(\form^n,\dom^n)$ is recurrent for each $n$ and converges in the sense of Mosco to the 
transient Dirichlet form $(\form, \dom)$.
\item[\sf (ii)] {\rm (transient ones to recurrent one)}
If we set 
$$
\alpha_n(u)= 1- \Bigl(\log \bigl(u+e^2\bigr)\Bigr)^{-(1-1/n)}, \quad  
\alpha(u) = 1- \Bigl(\log \bigl(u+e^2\bigr)\Bigr)^{-1}
$$ for $u\ge 0$ and $n \in \mathbb N$, then $(\form^n,\dom^n)$ is transient for each $n$ and converges in the sense of Mosco to the 
recurrent Dirichlet form $(\form, \dom)$.
\end{itemize}
\end{proposition}

\noindent
The organization of the paper is as follows. In the next section, we recall the Mosco convergence 
and give sufficient conditions for the Mosco convergence of the three types of Dirichlet forms. 
In Section \ref{sec: Instability}, we give several examples where global path properties are not preserved under the Mosco convergence. In Appendix, we give a necessary and sufficient condition for the recurrence 
of a class of symmetric stable type L\'evy processes.

\section{Mosco Convergence of Symmetric Dirichlet Forms on $L^2(\rd)$}

In the first part of this section, we briefly recall the notion of Mosco convergence
following \cite{M94}. 
For a closed form $(\form, \dom)$ on a Hilbert space ${\cal H}$, let $\form(u,u)=\infty$ for every 
$u\in {\cal H}\setminus\dom$. Here a closed form means a nonnegative definite symmetric closed form on 
${\cal H}$, not necessarily densely defined. 

\begin{definition} \label{defi: Mosco}{\rm
A sequence of closed forms $\form^n$ on a Hilbert space ${\cal H}$ is said to be convergent to $\form$ in 
the sense of Mosco if the following two conditions are satisfied:
\begin{itemize}
\item[{\sf (M1)}]  for every $u$ and every sequence  $\{u_n\}$ converging to $u$ {\it weakly} in ${\cal H}$, 
$$
\liminf_{n\to\infty} \form^n(u_n, u_n) \ge \form (u,u);
$$

\item[{\sf (M2)}]  for every $u$ there exists a sequence $\{u_n\}$ converging to $u$ in ${\cal H}$ 
so that 
$$
\limsup_{n\to\infty} \form^n(u_n, u_n) \le \form (u,u).
$$
\end{itemize}
}
\end{definition}

In \cite{M94}, Mosco showed that a sequence of closed forms $\form^n$ on ${\cal H}$ is converging to $\form$ 
in the sense of Mosco if and only if the resolvents associated with $\form^n$ converges to the resolvent associated 
with $\form$ {\it strongly} on ${\cal H}$, and hence the semigroups associated with $\form^n$ converges to 
the semigroup associated with $\form$ {\it strongly} on ${\cal H}$.

\subsection{Convergence of Symmetric Strong Local Dirichlet forms}

Consider a sequence of forms 
$$
\form^n(u,v)=\int_{\rd} \langle A_n(x) \nabla u(x), \nabla v(x)\rangle dx 
$$
for some functions $u,v$ in $L^2(\rd)$, where $A_n(x)=(a_{ij}^n(x))$ are $d\times d$ symmetric 
matrix valued functions satisfying {\bf Assumption (A)}.

\medskip
Under the assumption {\bf (A)}, the forms $(\form^n, C_0^{\infty}(\rd))$ for each $n$ and 
$(\form, C_0^{\infty}(\rd))$ are  Markovian closable forms on $L^2(\rd)$.
They become regular symmetric Dirichlet forms $(\form^n, \dom^n)$ and $(\form, \dom)$ on $L^2(\rd)$ 
(see \cite{FOT11}). Note that we set $\form^n(u,u)=\infty$ if $u \in L^2(\rd)\setminus \dom^n$.
We first give a simple lemma which is used in showing that $\form^n$ converges to $\form$ in the sense of Mosco.

\begin{lemma}
Suppose that Assumption {\sf (A)} holds. For any compact set $K\subset \rd$, 
 there exists a subsequence $\{n_k\}_k$ so that  $\int_K\|A_{n_{k}}^{1/2}(x)-A^{1/2}(x)\|^2dx 
 \to 0$ as $k \to\infty$.
\label{lem:diffu-m1-1}
\end{lemma} 

\noindent
{Proof:} \ Since $A_n(x)$ is a non-negative definite matrix for each $x$, there exists a nonnegative 
definite matrix $A_n^{1/2}(x)$ so that $\bigl(A_n^{1/2}(x)\bigr)^2=A(x)$. Then by the uniform boundedness 
of $A_n^{1/2}$ on the compact set $K$, we have 
\begin{align*}
\|A_n^{1/2}-A^{1/2}\|^2 \le \|A_n^{1/2}\|^2+2\|A_n^{1/2}\| \cdot \|A^{1/2}\| + \|A^{1/2}\|^2 < \infty.
\end{align*}
By {\sf(A2)}, there exists a subsequence $\{n_k\}_k$ so that $A_{n_k}(x) \to A(x)$ for {\it a.e.} $x \in K$ 
with respect to the matrix norm $\|\cdot\|$.  By general theory of linear operators, we can check that 
$A_n^{1/2}(x)\to A^{1/2}(x)$ in almost everywhere $x \in K$ with respect to $\|\cdot\|$. Thus, by using 
the dominated convergence theorem, we finish the proof. \ \hfill \fbox{}

\bigskip
We now prove Theorem \ref{thm: diffusion}:

\noindent
{Proof of Theorem \ref{thm: diffusion}}:  We first show {\sf (M1):}  Take $u\in L^2(\rd)$ and any 
sequence $\{u_n\}\subset L^2(\rd)$ so that $u_n$ converges to $u$ in $L^2$ weakly. We may assume that 
$\dis \liminf_{n\to\infty} \form^n(u_n, u_n)<\infty$. Taking a subsequence $\{n_k\}$, we have
\begin{align}
\infty> \liminf_{n\to\infty} \form^n(u_n, u_n) =\lim_{k\to\infty} \form^{n_k}(u_{n_k},u_{n_k})
&=\lim_{k\to\infty} \int_{\rd} |A_{n_k}^{1/2}\nabla u_{n_k}(x)|^2dx\label{eq:diffu-m1}.
\end{align}
Let us set $A_n^{1/2}(x)=(b_{ij}^n(x))$ and  $A^{1/2}(x)=(b_{ij}(x))$. By (\ref{eq:diffu-m1}), there exists $w_i\in L^2(\rd)$ so that, by taking a subsequence of 
$\{u_{n_k}\}$ if necessary,  $\sum_{j=1}^db_{ij}^{n_k}\partial_j u_{n_k}$ converges weakly to 
$w_i$ in $L^2(\rd)$ for each $i$.

We now  show that $w_i=\sum_{j=1}^db_{ij}\partial_j u$. To this end, take any 
$\eta \in C_0^{\infty}(\rd)$. Then we find that 
\begin{align*}
& \hspace*{-10pt} \int_{\rd} \bigl(w_i-\sum_{j=1}^db_{ij}(x)\partial_j u(x)\bigr)\eta(x)dx  \\
&  \ = \ \sum_{j=1}^d\int_{\rd} \bigl(w_i-b_{ij}^{n_k}(x)\partial_j u_{n_k}(x)\bigr)\eta(x)dx 
 + \ \sum_{j=1}^d\int_{\rd} \bigl(b_{ij}^{n_k}(x)\partial_j u_{n_k}(x)
-b_{ij}(x)\partial_j u_{n_k}(x)\bigr)\eta(x)dx  \\
& \quad \quad + \ \sum_{j=1}^d\int_{\rd} \bigl(b_{ij}(x)\partial_j u_{n_k}(x)-
b_{ij}(x)\partial_j u(x)\bigr)\eta(x)dx  \ =: \ {\rm (I)}_{k} +{\rm (II)}_{k} +
{\rm (III)}_{k}.
\end{align*}
We know that ${\rm (I)}_{k}$ converges to zero by definition. 
Now let us denote by $K$ the support of the function $\eta$. Then 
$$
{\rm (II)}_{k} \ = \ 
\sum_{j=1}^d\int_{K} \bigl(b_{ij}^{n_k}(x)
-b_{ij}(x)\bigr)\partial_j u_{n_k}(x)\eta(x)dx.
$$
By Lemma \ref{lem:diffu-m1-1}, taking a subsequence if necessary, $b_{ij}^{n_k}$ converges to $b_{ij}$ 
in $L^2(K)$.  Since $u_n$ converges weakly to $u$, $(\partial_j u_{n_k})$ also converges weakly to 
$(\partial_j u)$ in $L^2(\rd)$. Thus ${\rm (II)}_{k}$ converges to $0$ as $k\to\infty$ by the Schwartz 
inequality and $L^2$-boundedness of the weak-convergent sequence $\{\partial_j u_n\}_n$. The third term 
${\rm (III)}_{k}$ converges to $0$ since $\partial_i u_{n_k}$ converges weakly to $\partial_i u$ in 
$L^2(\rd)$ and $b_{ij}\eta \in L^2(\rd)$. Thus $w_i=\sum_{j=1}^db_{ij}\partial_j u$ holds for each 
$i=1,2,\ldots, d$.  Hence we have $\sum_{j=1}^db_{ij}^{n_k}\partial_j u_{n_k}$ converges weakly to 
$\sum_{j=1}^db_{ij}\partial_j u$, and we conclude that 
$$
\liminf_{n\to\infty} \form^n(u_n,u_n) =\lim_{k\to\infty} \int_{\rd} 
|A^{1/2}_{n_k}\nabla u_{n_k}(x)|^2dx \ge \int_{\rd} |A^{1/2}\nabla u(x)|^2dx=\form(u,u). 
$$
\bigskip
\noindent
We second show {\sf (M2):}  It is enough to show for $u\in \dom$. Since $C_0^{\infty}(\rd)$ is a (common) 
core for the Dirichlet forms $(\form^n,\dom^n)$, there exists a sequence $\{u_{\ell}\}\subset 
C_0^{\infty}(\rd)$ so that 
\begin{align} \label{conv: MD1}
\lim_{\ell \to\infty} \form_1(u_{\ell}-u,u_{\ell}-u) =
\lim_{\ell \to\infty} \Bigl( \int_{\rd} |A\nabla u_{\ell}(x)- A\nabla u(x)|^2dx + 
\int_{\rd} |u_{\ell}(x)-u(x)|^2 dx\Bigr) = 0.
\end{align}
Since any norms in the space of $d\times d$-matrices are equivalent, by Lemma \ref{lem:diffu-m1-1} and 
taking a subsequence if necessary, it follows that for each $\ell\in{\mathbb N}$, 
\begin{align*}
 \int_{\rd} \Bigl|A_n^{1/2}\nabla u_{\ell}(x)- 
A^{1/2} \nabla u_{\ell}(x)\Bigr|^2 dx 
 &\le \int_{\tilde{K}_{\ell}} \|A_n^{1/2}(x)-A^{1/2}(x)\|^2_{op}|\nabla u_{\ell}|^2_{\rd}(x)dx 
\\
&\le C\||\nabla u_{\ell}|_{\rd}\|^2_{\infty} \int_{\tilde{K}_{\ell}} \|A_n^{1/2}(x)-A^{1/2}(x)\|^2dx
\\
&\to 0 \quad {\rm as} \quad n\to\infty,
\end{align*}
where $C>0$ denotes some constant such that $\|\cdot\|_{op} \le C \|\cdot\|$ and $\|A\|_{op}$ means the operator norm of $A$: 
$\|A\|=\sup_{u: |u|_{\rd}\le 1}|Au|_{\rd}/|u|_{\rd}$. This gives us that 
$$
\lim_{n\to \infty} \form^n(u_{\ell},u_{\ell})=\form(u_{\ell},u_{\ell}), \quad \ell \in{\mathbb N}.
$$
Thus, with \eqref{conv: MD1}, we have
$$
\lim_{\ell\to \infty} \lim_{n\to \infty} \form^n(u_{\ell},u_{\ell})=\form(u,u).
$$
This shows {\sf (M2)}. \hfill \fbox{}


\subsection{Convergence of Symmetric Translation-Invariant Dirichlet Forms}

Let $\{\nu_t\}_{t>0}$ be a sequence of probability measures on $\rd$ of a continuous symmetric convolution 
semigroup:
$$
\left\{
\begin{array}{l}
\nu_t *\nu_s(A)=\nu_{t+s}(A), \quad t,s>0, \ A\in{\cal B}(\rd), \\
 \vspace{-10pt} \\
\nu_t(A)=\nu_t(-A), \quad A\in{\cal B}(\rd), \\
 \vspace{-10pt} \\
\nu_t \longrightarrow \delta \ \  {\rm weakly},
\end{array}
\right.
$$
where $\nu_t*\nu_s$ denotes the convolution of $\nu_t$ and $\nu_s$ \, $\bigl( \, \nu_t*\nu_s(A)
:=\int \nu_t(A-x)\nu_s(dx), \ A\in{\cal B}(\rd) \bigr)$ and $\delta$ is the Dirac measure concentrated at 
the origin. Define the kernels by 
$$
p(t,x, A):=p_t(x,A):= \nu_t(A-x), \quad t>0, \ x\in\rd, \ A\in{\cal B}(\rd),
$$
then $\{p_t(x,A) ; t>0, x\in \rd, \ A\in {\cal B}(\rd)\}$ is a Markovian transition function which is 
symmetric with respect to the Lebesgue measure in the following sense:
$$
\int_{\rd} p_tf(x)g(x)dx=\int_{\rd} f(x)p_tg(x)dx, \quad f, g\in {\cal B}(\rd)_+.
$$
According to the L\'evy-Khinchin formula, we see that a continuous symmetric convolution 
semigroup $\{\nu_t, t>0\}$ is characterized by a pair $(S, n)$ satisfying (\ref{LK-1}) and 
(\ref{LK-2}) through the formula (\ref{LK-formula}).

Now let $\{\varphi_n\}$ be a sequence of the characteristic functions defined by symmetric convolution 
semigroups $\{ \nu_t^n, t>0\}_{n\in\mathbb N}$: 
$$
e^{-t \varphi_n(x)} := \widehat{\nu}_t^n(x) 
\Bigl(= \int_{\rd}e^{i\langle x,y\rangle} \nu_t^n(dy)\Bigr), \quad x\in \rd.
$$
Let $\varphi$ be also a characteristic function defined by a symmetric convolution 
semigroup $\{\nu_t, t>0\}$. The Dirichlet forms corresponding $\nu^n_t$ are defined by 
$$
\left\{
\begin{array}{rcl}
\form^n(u,v) &=& \dis \int_{\rd} \widehat{u}(x)\overline{\widehat{v}}(x) \varphi_n(x)dx \\ 
& & \vspace*{-6pt} \\ 
\dom^n  &=& \dis 
\Bigl\{ u \in L^2(\rd) : \ \int_{\rd} \bigl|\widehat{u}(x)\bigr|^2 \varphi_n(x)dx<\infty
\Bigr\}.
\end{array}
\right.
$$
We set that for each $n$, $\form^n(u,u)=\infty$ if $u \in L^2(\rd)\setminus \dom^n$. 
We assume Assumption {\bf (B)}.

\bigskip
Then we show Theorem \ref{Mosco:Levy}:

\noindent
{Proof of Theorem \ref{Mosco:Levy}:} We first show {\sf (M1):} \ Take any $u\in L^2(\rd)$ and any 
sequence $\{u_n\} \subset L^2(\rd)$ for which $u_n$ converges to $u$ weakly in $L^2(\rd)$. 
We may assume $\dis \liminf_{n\to\infty} \form^n(u_n,u_n) <\infty$. 

According to the Parseval formula, note that $u_n$ converges to $u$ weakly in $L^2$ if and only if 
$\widehat{u}_n$ converges to $
\widehat{u}$ wealky in $L^2$. Here $\widehat{u}$ denotes the Fourier transform of $u$. Thus 
$$
\infty \ > \ \liminf_{n\to\infty} \form^n(u_n,u_n) 
\ = \ \liminf_{n\to\infty} \int_{\rd} \bigl|\widehat{u}_n(x)\bigr|^2 \varphi_n(x)dx
$$
implies that there exist a subsequence $\{n_k\}$ and an element $w\in L^2(\rd)$ 
so that $\widehat{u}_{n_k} \cdot \sqrt{\varphi_{n_k}}$ converges to $w$ weakly in $L^2(\rd)$.
We now show that $w=\widehat{u}\cdot \sqrt{\varphi}$.  For any $v \in C_0^{\infty}(\rd)$,  we see
that 
\begin{align*}
 & \Bigl| \int_{\rd} \bigl(w(x)-\widehat{u}(x)\sqrt{\varphi(x)}\bigr) v(x)dx \Bigr| \\
&  \ \le \ \Bigl| \int_{\rd} \bigl(w(x)-\widehat{u}_{n_k}(x)\sqrt{\varphi_{n_k}(x)}\bigr) v(x)dx \Bigr| 
+ \ \Bigl| \int_{\rd} \bigl(\widehat{u}_{n_k}(x)\sqrt{\varphi_{n_k}(x)}-
\widehat{u}(x) \sqrt{\varphi_{n_k}(x)} \bigr) v(x)dx \Bigr|  \\
&  \quad  \quad + \ 
 \ \Bigl| \int_{\rd} \bigl(\widehat{u}(x) \sqrt{\varphi_{n_k}(x)} -\widehat{u}(x) \sqrt{\varphi(x)}\bigr) 
 v(x)dx \Bigr| \\
& \ =: \ {\rm (I)}_k + {\rm (II)}_k +{\rm (III)}_k. 
\end{align*}
The first term ${\rm (I)}_k$ converges to $0$ by definition.  For the second term ${\rm (II)}_k$, using 
the condition {\bf  (B)} and the inequality \textcolor{black}{$|\sqrt{a}-\sqrt{b}| \le \sqrt{|a-b|}$}, 
we have that $\sqrt{\varphi_{n}}v$ converges to $\sqrt{\varphi}v$ in $L^2(\rd)$. Thus the second term 
${\rm (II)}_k$ converges to zero by the Schwartz inequality and $L^2$-boundedness of 
$\{\widehat{u}_{n_k}\}_k$. For the third term ${\rm (III)}_k$,
\begin{align*}
{\rm (III)}_k & \ = \ 
\Bigl| \int_{\rd} \bigl(\sqrt{\varphi_{n_k}(x)} -\sqrt{\varphi(x)}\bigr) \widehat{u}(x)v(x)dx \Bigr|  \\
& \ \le \ \|\widehat{u}\|_{L^2}
\sqrt{ \int_{\rd} \Bigl( \sqrt{\varphi_{n_k}(x)}v(x) -\sqrt{\varphi(x)}v(x) \ \Bigr)^2 dx } 
\ \longrightarrow 0 \ \ {\rm as}
 \ n\to \infty.
\end{align*}
Thus we conclude that $w=\widehat{u}\sqrt{\varphi}$. Hence we find that 
$$
\liminf_{n\to \infty} \form^n(u_n,u_n) = \lim_{k\to \infty}
\form^{n_k}(u_{n_k}, u_{n_k}) = \lim_{k\to \infty} \int_{\rd} 
\bigl|\widehat{u}_{n_k}\bigr|^2 \varphi_{n_k}dx \ge  \int_{\rd} 
\bigl|\widehat{u}\bigr|^2 \varphi dx = \form(u,u).
$$
This shows (M1).

\bigskip
\noindent
We second show {\sf (M2):}  It is enough to show for $u \in \dom$. Since $C_0^{\infty}(\rd)$ is a (common) 
core for the Dirichlet forms, there exists a sequence $\{\tilde{u}_n\}$ of $C_0^{\infty}(\rd)$ such that 
\begin{equation}
\lim_{n\to\infty}\form_1(\tilde{u}_n-u,\tilde{u}_n-u)=
\lim_{n\to\infty} \Bigl(\int_{\rd} \bigl|\widehat{\tilde{u}}_n(x)-\widehat{u}(x)\bigr|^2 
\varphi(x)dx + \int_{\rd} \bigl(\tilde{u}(x)-u(x)\bigr)^2dx \Bigr)=0.
\label{approx-1}
\end{equation}
We now take a sequence $\{\chi_n\}$ of $C_0^{\infty}(\rd)$ satisfying 
$$
\chi_n(x)=\chi_n(-x), \quad  0\le \chi_n(x) \le \chi_{n+1}(x) \le 1, \ \ 
n \in{\mathbb N}, \ \  \lim_{n\to \infty} \chi_n(x)=1, \ x\in\rd.
$$
For any $n,\ell \in{\mathbb N}$, \ set $u_{n, \ell}(x)=\tilde{u}_n*\check{\chi}_{\ell}(x) = 
\int_{\rd}\check{\chi}_{\ell}(x-y)\tilde{u}_n(y)dy, \ x\in\rd$. Here the inverse Fourier transform of 
$\chi_{\ell}$ is denoted by $\check{\chi_{\ell}}$. Since
$$
\| \tilde{u}_n * \check{\chi}_{\ell} - u *\check{\chi}_{\ell} \|_{L^2} 
\ = \ \| ( \widehat{\tilde{u}}_n  - \widehat{u} ) \cdot \chi_{\ell} \|_{L^2} \ 
\le \ \| \widehat{\tilde{u}}_n  - \widehat{u} \|_{L^2} \ = \ \| \tilde{u}_n -u\|_{L^2} \ \ 
\longrightarrow \ 0 \quad {\rm as} \ n\to \infty
$$
for each $\ell$ and 
$$
\|u*\check{\chi}_{\ell}-u\|_{L^2}=
\|\widehat{u}\cdot \chi_{\ell}-\widehat{u}\|_{L^2} \ \ \longrightarrow \  0 
\quad {\rm as} \ \ell \to \infty,
$$
then we have that 
$$
\lim_{\ell\to \infty} \lim_{n\to \infty} \|u_{n, \ell}-u\|_{L^2}=0.
$$
On the other hand, we see that from (\ref{approx-1}), 
$\widehat{\tilde{u}}_{n} \cdot \chi_{\ell} \sqrt{\varphi_n}
=\widehat{\tilde{u}_n*\check{\chi}_{\ell}}\sqrt{\varphi_n}$ converges to 
$\widehat{u} \cdot \chi_{\ell} \sqrt{\varphi}
=\widehat{u*\check{\chi}_{\ell}}\sqrt{\varphi}$ 
in $L^2(\rd)$ for any $\ell$. 
In fact, using the inequalities $(a-b)^2 \le 2a^2+2b^2$, \textcolor{black}{$|\sqrt{a}-\sqrt{b}| \le 
\sqrt{|a-b|}$} and the condition {\bf  (B)}, we have
\begin{align*}
 & \hspace*{-1cm} \int_{\rd} \Bigl(\widehat{\tilde{u}}_n\cdot \chi_{\ell} \sqrt{\varphi_n} 
-\widehat{u}\cdot \chi_{\ell}\sqrt{\varphi} \Bigr)^2 dx  \\
& \ \le \ 2\int_{\rd} \Bigl(\widehat{\tilde{u}}_n\cdot \chi_{\ell}\sqrt{\varphi_n} 
-\widehat{\tilde{u}}_n\cdot \chi_{\ell}\sqrt{\varphi} \Bigr)^2 dx + 2
\int_{\rd} \Bigl(\widehat{\tilde{u}}_n\cdot \chi_{\ell}\sqrt{\varphi} 
-\widehat{u}\cdot \chi_{\ell}\sqrt{\varphi} \Bigr)^2 dx \\ 
& \ \le \ 2\int_{\rd} \widehat{\tilde{u}}_n^2\chi_{\ell}^2|\varphi_n 
-\varphi | dx  + \textcolor{black}{2
\int_{\rd} \Bigl(\widehat{\tilde{u}}_n\cdot \chi_{\ell}\sqrt{\varphi} 
-\widehat{u}\cdot \chi_{\ell}\sqrt{\varphi} \Bigr)^2 dx} \\
& \longrightarrow \ 0 \quad {\rm as} \ n\to\infty.
\end{align*}
Thus we find that 
$$
\lim_{n\to \infty} \form^n(u_{n,\ell}, u_{n,\ell})=\lim_{n\to \infty}
\int_{\rd} \bigl|\widehat{\tilde{u}_n*\check{\chi}_{\ell}}\bigr|^2 \varphi_ndx
=\int_{\rd} \bigl|\widehat{u*\check{\chi}_{\ell}}\bigr|^2 \varphi dx
=\int_{\rd} \bigl|\widehat{u}\bigr|^2 \chi_{\ell}^2 \varphi dx
$$
and 
$$
\lim_{\ell \to \infty} 
\int_{\rd} \bigl|\widehat{u}\bigr|^2 \chi_{\ell}^2 \varphi dx
=\int_{\rd} \bigl|\widehat{u}\bigr|^2\varphi dx
=\form(u,u).
$$
These imply that 
$$
\lim_{\ell \to \infty} 
\lim_{n\to \infty} \form^n(u_{n,\ell}, u_{n,\ell})
=\form(u,u)
\quad {\rm and}  \quad 
\lim_{\ell \to \infty} 
\lim_{n\to \infty} \int_{\rd} \bigl| u_{n,\ell} -u\bigr|^2 dx =0.
$$
Therefore, by the diagonalization argument, we can find a sequence $\{\ell(n)\}_n$ so that 
$$
\ell(n)<\ell(n+1) \nearrow \infty \ (n\to\infty), \quad 
\lim_{n\to\infty} \form^{n} (u_{n,\ell(n)},u_{n,\ell(n)})=\form(u,u)
$$
and then {\sf (M2)} is shown. \hfill \fbox{}

\subsection{Convergence of Symmetric Jump-Type Dirichlet forms}
Let $\wt{J}(x,y)$ be a non-negative symmetric function on $\rd\times \rd \setminus {\sf diag}$ 
satisfying 
\begin{align} \label{eq: Loc}
x\mapsto \int_{y\not=x} \Bigl(1\wedge d(x,y)^2\Bigr) \wt{J}(x,y)dy \in L^1_{\sf loc}(\rd).
\end{align} 
Consider the following quadratic form $\wt\form$ on $L^2(\rd)$: 
$$
\wt{\form}(u,v) = \frac{1}{2}\iint_{x\not=y} \bigl(u(x)-u(y)\bigr)
\bigl(v(x)-v(y)\bigr)\wt{J}(x,y)dxdy
$$
for some functions $u,v\in L^2(\rd)$. 
Under the condition \eqref{eq: Loc}, it is known that $(\form, C_0^{\sf lip}(\rd))$ 
is a closable Markovian symmetric form on $L^2(\rd)$.
Thus taking the closure of $C_0^{\sf lip}(\rd)$ with respect to 
$\sqrt{\wt\form_1}$, we find that $(\wt\form,\wt\dom)$ is a regular Dirichlet form.  

Now take $J_n(x,y)$ and $J(x,y)$  non-negative symmetric functions on
 $\rd\times \rd \setminus {\sf diag}$ satisfying \eqref{eq: Loc} in place of $\wt J(x,y)$ and then 
 consider regular symmetric Dirichlet forms $\form^n$ and $\form$ of pure jump type on 
 $L^2(\rd)$ as follows:
$$
\left\{
\begin{array}{rcl}
\form^n(u,v) &=& \dis{\frac{1}{2}\iint_{x\not=y} \bigl(u(x)-u(y)\bigr)
\bigl(v(x)-v(y)\bigr)J_n(x,y)dxdy, } \\
& & \vspace{-10pt} \\
\dom^n &=&  \overline{C_0^{\sf lip}(\rd)}^{\sqrt{\form_1^n}},
\end{array}
\right.$$
and 
$$
\left\{
\begin{array}{rcl}
\form(u,v) &=& \dis \frac{1}{2}\iint_{x\not=y} \bigl(u(x)-u(y)\bigr)
\bigl(v(x)-v(y)\bigr)J(x,y)dxdy, \\
& & \vspace{-10pt} \\
\dom &=&  \overline{C_0^{\sf lip}(\rd)}^{\sqrt{\form_1}}.
\end{array}
\right.
$$

\medskip
We assume Assumption {\bf (C)}. 
Then we prove Theorem \ref{thm: JtoJ}.

\bigskip
\noindent
{Proof of Theorem \ref{thm: JtoJ}:} \ We have to check the following two conditions:
\begin{itemize}
\item[{\sf (M1)}] \ for any $u\in L^2(\rd)$ and $\{u_n\}\subset L^2(\rd)$ which 
converges to $u$ weakly in $L^2(\rd)$,
$$
\liminf_{n\to\infty} \form^n(u_n,u_n)\ge \form(u,u).
$$
\item[{\sf (M2)}] \ for any $u\in L^2(\rd)$, there exists a sequence 
$\{u_n\}\subset L^2(\rd)$ which converges to $u$ in $L^2(\rd)$ such that 
$$
\limsup_{n\to\infty} \form^n(u_n,u_n) \le \form(u,u).
$$
\end{itemize}

\noindent
{Proof of {\sf (M1)}:} \ Suppose that 
\begin{itemize}
\item[(1)] \ $u_n$ is weakly convergent to $u$ in $L^2(\rd)$ and 
\item[(2)] \ $\dis{\liminf_{n\to\infty} \iint_{x\not=y} 
\bigl(u_n(x)-u_n(y)\bigr)^2 J_n(x,y)dxdy<\infty}$.
\end{itemize} 

We may assume that 
$\dis{\lim_{n\to\infty} \iint_{x\not=y} 
\bigl(u_n(x)-u_n(y)\bigr)^2 J_n(x,y)dxdy<\infty}$.

Then for each $n$, put ${\bar u}_n(x,y)=\bigl(u_n(x)-u_n(y)\bigr)\sqrt{J_n(x,y)}$ 
\ for $(x,y) \in \rd\times \rd \setminus {\sf diag}$. 
Then $\{{\bar u}_n\}$ are bounded sequence in $L^2(\rd{\times}\rd \setminus{\sf diag}; 
dx{\otimes}dy)$, and so there exists a subsequence $\{{\bar u}_{n_k}\}$ 
which converges to some element ${\bar u}$ {\it weakly} in 
$L^2(\rd{\times}\rd\setminus{\sf diag}; 
dx{\otimes}dy)$. 

We now claim that 
$$
{\bar u}(x,y)=\bigl(u(x)-u(y)\bigr) \sqrt{J(x,y)}, \quad 
dx{\otimes}dy\mbox{-a.e.} \ (x,y) \ {\rm with} \ x\not=y.
$$ 
For any nonnegative $v \in C_0(\rd \times \rd \setminus {\sf diag})$ and for any $n_k$, we see

\begin{align*}
 & \hspace*{-10pt} \Bigl\vert \iint_{x\not=y} \Bigl({\bar u}(x,y)-
\bigl(u(x)-u(y)\bigr)\sqrt{J(x,y)}\Bigr)v(x,y) dxdy \Bigr\vert \\
 &  \le \ \Bigl\vert \iint_{x\not=y} \Bigl({\bar u}(x,y)-
\bigl(u_{n_k}(x)-u_{n_k}(y)\bigr)\sqrt{J_{n_k}(x,y)} \, \Bigr)v(x,y) dxdy \Bigr\vert  \\
& \quad +  \Bigl\vert \iint_{x\not=y} \bigl(u_{n_k}(x)-u_{n_k}(y)\bigr)
\bigl(\sqrt{\mathstrut J_{n_k}(x,y)}-\sqrt{\mathstrut J(x,y)} \Bigr) v(x,y) dxdy
\Bigr\vert  \\
& \quad +  \Bigl\vert \iint_{x\not=y} 
\Bigl( \bigl(u_{n_k}(x)-u_{n_k}(y)\bigr) -\bigl(u(x)-u(y)\bigr) \Bigr)
\sqrt{\mathstrut J(x,y)} v(x,y) dxdy \Bigr\vert \\
& =: {\rm (I)}_{n_k} + {\rm (II)}_{n_k} + {\rm (III)}_{n_k}. 
\end{align*}
Since ${\bar u}_n$ converges to ${\bar u}$ {\it weakly} in 
$L^2(\rd{\times}\rd \setminus {\sf diag};dx\otimes dy)$, we see $\lim_{k\to\infty} {\rm (I)}_{n_k}=0$.  
By making use of the Schwarz inequality and Assumption {\bf (C)} and 
noting $\{u_{n_k}\}$ is a bounded sequence in $L^2(\rd;dx)$, we see 
\begin{align*}
{\rm (II)}_{n_k} \ & \le \ 
\sqrt{\iint_{x\not=y}  \hspace*{-5pt} \bigl(u_{n_k}(x){-}u_{n_k}(y)\bigr)^2
v(x,y) dxdy} \\
& \qquad \times  \sqrt{ \iint_{x\not=y} \hspace*{-5pt} 
\Bigl( \sqrt{J_{n_k}(x,y)}- \sqrt{J(x,y)}\,  \Bigr)^2 v(x,y)dxdy} \\
& \le \  \dis \| v\|_{\infty} \|u_{n_k}\|_{L^2}  \sqrt{2 \left(
\left\|\int_{\{x: (x,\cdot) \in {\rm supp}[v]\}} \hspace*{-10pt} |v(x, \cdot)|dx\right\|_{\infty}  +
\left\|\int_{\{y: (\cdot, y) \in {\rm supp}[v]\}} \hspace*{-10pt} |v(\cdot, y)|dy\right\|_{\infty}
\right)}  \\
& \qquad \times 
\sqrt{\iint_{ {\rm supp}[v]} |J_{n_k}(x,y)- J(x,y)|  dxdy}   \\
& \longrightarrow \  0  \quad {\rm as} \ n_k\to\infty.
\end{align*}
Here we uesd elementary inequalities in the second inequality above: $(a-b)^2 \le 2(a^2+b^2)$ and
$|\sqrt{a}-\sqrt{b}| \le \sqrt{|a-b|}$ for $a, b\ge 0$. As for ${\rm (III)}_{n_k}$, 
note that both
$$
\varphi(x)=\int_{y\not=x} \sqrt{J(x,y)} v(x,y)dy, \  x \in \rd, \quad  
\psi(y)=\int_{y\not=x}\sqrt{J(x,y)}v(x,y)dx, \  y\in \rd
$$ 
are in $L^2(\rd)$. So we see 
$$
{\rm (III)}_{n_k} \le \Bigl\vert 
\int_{\rd} \bigl(u_{n_k}(x)-u(x)\bigr) \varphi(x)dx\Bigr\vert + 
\Bigl\vert 
\int_{\rd} \bigl(u_{n_k}(y)-u(y)\bigr) \psi(y)dy \Bigr\vert
$$
goes to $0$ when $n_k\to\infty$. Thus we see
$$
{\bar u}(x,y) =\bigl(u(x)-u(y)\bigr)\sqrt{J(x,y)} \quad dx\otimes dy\mbox{-a.e.} \ \  
(x,y)  \  {\rm with} \ x\not=y.
$$
Hence 
$$
\liminf_{n\to\infty} \form^n(u_n,u_n) 
\ge \iint_{x\not=y} \bigl(u(x)-u(y)\bigr)J(x,y)dxdy
=\form(u,u).
$$ 
\ \rightline{\fbox{}}

\bigskip
\noindent
{Proof of {\sf (M2):}} \  Since $C_0^{\sf lip}(\rd)$ is a common core for $(\form^n, \dom^n)$, 
it is enough to show {\sf (M2)} for functions in $C_0^{\sf lip}(\rd)$. 
Take any $u\in C_0^{\sf lip}(\rd)$. Put $u_n=u$ for each $n$, then 
$u_n$ converges to $u$ in $L^2(\rd)$. Denote by $K$ the support of $u$ and 
take compact set $F$ with $K\subset F$ and 
\begin{align}
d(K,F^c)=\inf\{d(x,y) : \ x\in K, \ y \in F^c\}\ge 1. \label{eq: SET}
\end{align}
Then 
$$
\begin{array}{rcl}
\form^n(u_n,u_n)=\form^n(u,u)
 &=& \dis{\iint_{x\not=y} \bigl(u(x)-u(y)\bigr)^2J_n(x,y)dxdy } \\
& & \vspace{-8pt} \\
& = &  \dis{ \iint_{F\times F \setminus {\sf diag}} \underline{\hspace{1cm}} \, 
 dxdy + 2\iint_{K \times F^c} \underline{\hspace{1cm}}\,dxdy } \\
& & \vspace{-8pt} \\
& \equiv & {\rm (I)}_n + 2{\rm (II)}_n. \\
\end{array}
$$
We first estimate $(I)_n$. \ For all $n\in{\mathbb N}$ and $dx{\otimes}dy$-a.a. $(x,y) \in K\times K \setminus {\sf diag}$, 
we see that integrand in ${\rm (I)}_n$ is  $(u(x)-u(y))^2 J_n(x,y)$ and is 
bounded by $M \bigl(1\wedge d(x,y)^2\bigr) \wt J(x,y)$ from above, where $M$ is the 
maximum of the Lipschitz constant of $u$ and $\|u\|_{\infty}$. 
Then the fact that the function $\bigl(1\wedge d(x,y)^2\bigr) \tilde{J}(x,y)$ is 
integrable on the set $F \times F \setminus {\sf diag}$ and Assumption {\bf (C)} imply that 
$$
\lim_{n\to\infty} {\rm (I)}_n= \iint_{F\times F \setminus {\sf diag}} 
\bigl(u(x)-u(y)\bigr)^2 J(x,y)dxdy
$$
We next estimate $(II)_n$.  The integral $(II)_n$ is the following:
$$
\iint_{K\times F^c} u(x)^2 J_n(x,y)dxdy.
$$
For $dx{\otimes}dy$-a.a. $(x,y) \in K\times F^c$, we see from Assumption {\bf (C)}, 
we see
$$
u(x)^2 J_n(x,y) \le ||u||_{\infty}^2 \wt J(x,y)
$$
and the right hand side is integrable on the set $K\times F^c$ because of \eqref{eq: Loc} and  \eqref{eq: SET}.
So, as in the case ${\rm (I)}_n$, we have 
$$
\lim_{n\to\infty} {\rm (II)}_n=\iint_{K\times F^c} u(x)^2 
J(x,y)dxdy.
$$
Combining these two estimates, we have 
$$
\begin{array}{lcl}
\dis{\lim_{n\to\infty} \form^n(u,u)} &=& \dis{
\iint_{F \times F \setminus {\sf diag}} 
\bigl(u(x)-u(y)\bigr)^2 J(x,y)dxdy } \\
& & \vspace{-8pt} \\
& & \quad + \dis{  2 \iint_{K\times F^c} u(x)^2 
J(x,y) dxdy } \\
& & \vspace{-8pt} \\
&=&  \dis{\iint_{x\not=y} \bigl(u(x)-u(y)\bigr)^2 J(x,y)dxdy
= \form(u,u). } \\
\end{array}
$$
This means that the condition {\sf (M2)} holds. \quad \fbox{}

\section{Instability of Global Path Properties} \label{sec: Instability}
\subsection{Proof of Proposition \ref{prop: Cons/Exp}}
\bigskip
\noindent
{Proof of Proposition \ref{prop: Cons/Exp}:} 
By \cite[Theorem 2.2]{T89} and the Feller's test in \cite{Mc69}, in the case of (i), 
$(\form^n,\dom^n)$ is explosive and $(\form,\dom)$ is conservative, and,  in the case of (ii),  
$(\form^n,\dom^n)$ is conservative and $(\form,\dom)$ is explosive (see also, e.g., 
\cite[page 300]{FOT11}).  By Theorem \ref{thm: diffusion}, $(\form^n,\dom^n)$ converges to 
$(\form,\dom)$ in the sense of Mosco in the both cases (i) and (ii) and we finish the proof. \hfill \fbox{}

\subsection{Proof of Proposition \ref{prop: Rec/Tra}}

We first show the following lemma which is a sufficient condition for local 
$L^1$-convergence of $\varphi_n$ to $\varphi$:

\begin{lemma} \label{lem:Levyconv1}
If $\alpha_n(t) \to \alpha(t)$ for every $t \in [0,\infty)$, 
then $\varphi_n \to \varphi$ locally in $L^1(\rd)$.
\end{lemma}

\noindent
{Proof:} We show that, for any compact set $K \subset \rd$, $\int_K|\varphi_n(x)-\varphi(x)|dx 
\to 0$. Since $\varphi_n$ and $\varphi$ are continuous functions, 
$|\varphi_n(x)-\varphi(x)|$ is bounded on $K$ and, for making use of the dominated convergence theorem,  it suffices to show that 
$\varphi_n(x) \to \varphi(x)$ for a.e. $x \in K$. We see that 
\begin{align*}
   |\varphi_n(x)-\varphi(x)| 
    & \le \int_{\rd} \Bigl|1-\cos (\langle x, \xi\rangle)\Bigr|\cdot 
    \Bigl||\xi|^{-1-\alpha_n(|\xi|)}-|\xi|^{-1-\alpha(|\xi|)}\Bigr|d\xi \\
    & \le \int_{\rd} |1-\cos (\langle x, \xi\rangle)||\xi|^{-1-\alpha_n(|\xi|)}d\xi 
     + \int_{\rd} |1-\cos (\langle x, \xi\rangle)||\xi|^{-1-\alpha(|\xi|)}d\xi \\
    & \le \int_{\rd} |1-\cos (\langle x, \xi\rangle)|\max\{|\xi|^{-1-\overline{\alpha}}, 
       |\xi|^{-1-\underline{\alpha}}\}d\xi+ 
       \int_{\rd} |1-\cos (\langle x, \xi\rangle)|\max\{|\xi|^{-1-\overline{\alpha}},|\xi|^{-1-\underline{\alpha}}\}d\xi \\
    & < \infty.
\end{align*}
Since $\alpha_n(t) \to \alpha(t)$ as $n \to 0$ for any $t\in [0,\infty)$, 
it follows that $(1-\cos (\langle x, \xi\rangle))(|\xi|^{-1-\alpha_n(|\xi|)}-|\xi|^{-1-\alpha(|\xi|)}) \to 0$ as 
$n \to \infty$.  By the dominated convergence theorem, we see that $\varphi_n(x) \to \varphi(x)$ for 
a.e. $x \in \rd$. The proof is completed. \hfill \fbox{}

Now we show Proposition \ref{prop: Rec/Tra}: 

\bigskip
\noindent
{Proof of Proposition \ref{prop: Rec/Tra}:}  \ (i): By Theorem \ref{thm: U04}, we can verify that $(\form^n,\dom^n)$ is recurrent 
for any $n$ and $(\form, \dom)$ is transient. By Lemma \ref{lem:Levyconv1}, we have that 
$(\form^n,\dom^n)$ converges to $(\form, \dom)$ in the sense of Mosco.  (ii): The transience of 
$(\form^n,\dom^n)$ and the recurrence of $(\form, \dom)$ follow directly from Theorem \ref{thm: Rec}.
By Lemma \ref{lem:Levyconv1}, we have that $(\form^n,\dom^n)$ converges to $(\form, \dom)$ 
in the sense of Mosco and we finish the proof.
\hfill \fbox{}

\subsection{Proof of Proposition \ref{prop: Rec/Tra_J}}
\bigskip
\noindent
{Proof of Proposition \ref{prop: Rec/Tra_J}:} We use Proposition \ref{prop: Rec/Tra} and the comparison 
theorems of Dirichlet forms \cite[Theorem 1.6.4]{FOT11}.
\hfill \fbox{}

\section{Appendix: Sharpness of Recurrence Criteria for Symmetric L\'evy Processes}

It is well-known that a translation invariant symmetric stable process with an index $\alpha \ 
(0<\alpha \le 2)$ is recurrent if and only if $d=1\le \alpha \le 2$ or $d=\alpha=2$. 
The L\'evy kernel is given by $n(dh)=c |h|^{-d- \alpha}dh$ for some constant $c=c(d,\alpha)$ 
if $0<\alpha<2$. 

In this appendix,  we give a recurrent criteria for a calss of stable type L\'evy processes 
having the L\'evy measure $n(dh)=|h|^{-d-\alpha(|h|)} dh$, where $\alpha$ is a measurable function 
defined on $[0,\infty)$.  When $\alpha$ is a constant 
between $0$ to $2$, then this corresponds nothing but to a symmetric $\alpha$ stable process.
Consider also the following quadratic form:
\begin{align*}
\form (u,v) & \ = \ \iint_{h\not=0} \bigl(u(x+h)-u(x)\bigr)
\bigl(v(x+h)-v(x)\bigr) n(dh) dx \\
 & \ = \ \iint_{x\not=y}  \frac{\bigl(u(y)-u(x)\bigr)
\bigl(v(y)-v(x)\bigr)}{|x-y|^{d+\alpha(|x-y|)}} dxdy, \\
{\cal D}[\form] &  \ = \ \Bigl\{ u \in L^2(\rd) : \form(u,u)<\infty \Bigr\}.
\end{align*}
Then it is known that $(\form, {\cal D}[\form])$ is a symmetric, 
translation invariant Dirichlet form on $L^2(\rd)$ under the following 
condition:
$$ 
\int_{h\not=0} \bigl(1\wedge |h|^2\bigr) n(dh) = c_d \int_0^{\infty}
\bigl(1\wedge u^2 \bigr) u^{-1-\alpha(u)}du<\infty.
$$

\bigskip
In \cite{U04} (see also \cite{MUJ12, OU14}), we have shown the following theorem:

\begin{theorem}  \label{thm: U04}{\sf ({\it c.f.} Theorem 3.3 in \cite{U04}).} \ 
If the conditions 
$$
\limsup_{R\to \infty} R^{-2+d} \int_0^R u^{1-\alpha(u)}du<\infty
$$
and 
$$
\limsup_{R\to \infty} R^d \int_R^{\infty} u^{-1-\alpha(u)}du<\infty
$$
hold, then the process is recurrent.
\end{theorem}

In the case where $d=1$, we can show the following. Let $\varepsilon>0$ and set  
$$
\alpha(u) 
= 1- \Bigl(\log \bigl(u+e^2\bigr)\Bigr)^{-\vareps}, \quad u\ge 0
%
$$
for instance. Let us also consider the corresponding form:
\begin{align*}
\form (u,v) & \ = \ \iint_{h\not=0} \bigl(u(x+h)-u(x)\bigr)
\bigl(v(x+h)-v(x)\bigr) n(dh) dx \\
 & \ = \ \iint_{x\not=y}  \frac{\bigl(u(y)-u(x)\bigr)
\bigl(v(y)-v(x)\bigr)}{|x-y|^{1+\alpha(|x-y|)}} dxdy, \\
{\cal D}[\form] &  \ = \ \Bigl\{ u \in L^2({\mathbb R}) : \form(u,u)<\infty \Bigr\}.
\end{align*}

Then we have the following criterion for the recurrence:
\begin{theorem}  \label{thm: Rec}
The form/process is recurrent if and only if $\varepsilon \ge 1$.
\end{theorem}

\noindent
{Proof:} Though we have shown in \cite{U04} ({\it c.f} \cite{U02}) that 
the form is recurrent if $\varepsilon\ge 1$, we give the proof of it for 
reader's convenience. Namely, we estimate two integrals 
in the previous theorem in the case $d=1$.

For $R>e$, 
$$
R^{-1} \int_0^R u^{1-\alpha(u)}du \ = \  R^{-1} \int_0^R u^{(\log (u+e^2))^{-\varepsilon}}du.
$$
Since $\vareps \ge 1$, we find that 
\begin{align*}
R^{-1} \int_0^R u^{(\log (u+e^2))^{-\varepsilon}}du 
& \ = R^{-1} \int_0^{\sqrt{R}-e^2}u^{(\log (u+e^2))^{-\varepsilon}}du 
+ R^{-1} \int^R_{\sqrt{R}-e^2}u^{(\log (u+e^2))^{-\varepsilon}}du  \\
& \ \le \ R^{-1} \int_0^{\sqrt{R}} u^{2^{-\vareps}}du + R^{-1}
 \int^R_{\sqrt{R}-e^2}u^{(\frac 12 \log R)^{-\varepsilon}}du  \\
& \ \le \ \frac{R^{2^{-\vareps-1}-1}}{2^{-\vareps}+1} + 
R^{-1+ (\frac 12 \log R)^{-\varepsilon}} (R-\sqrt{R}+e^2) \\
& \ = \ \frac{2^{\vareps}}{(1+2^{\vareps})R^{1-2^{-\vareps-1}}}  + 
R^{(\frac 12 \log R)^{-\varepsilon}} \Bigl(1-\frac 1{\sqrt{R}}+\frac {e^2}{R}\Bigr).
\end{align*}
Since $\vareps \ge 1$, it follows that
$$
\log  R^{(\frac 12 \log R)^{-\varepsilon}} = 
(\frac 12 \log R)^{-\varepsilon} \cdot \log R =2^{\vareps} (\log R)^{1-\vareps} 
\to \left\{ {0  \quad {\rm if} \ \vareps>1 \atop 2^{\vareps} \quad {\rm if} \ \vareps=1} \right. 
\quad {\rm as} \ R\to\infty. 
$$
Thus we find that 
$$
\limsup_{R\to \infty} R^{-1} \int_0^R u^{1-\alpha(u)}du <\infty.
$$
We now estimate the second condition: For $R>\sqrt{e}$, 
\begin{align*}
R \int_R^{\infty} u^{-1-\alpha(u)}du & \ =  \ 
 R \int_R^{\infty} u^{-2+(\log (u+e^2))^{-\vareps}}du \ \le \ 
 R \int_R^{\infty} u^{-2+(2\log R)^{-\vareps}}du  \\
 & \ = \ R\Bigl[ \frac 1{-1 +(2\log R)^{-\vareps}} u^{-1+(2\log R)^{-\vareps}}\Bigr]_R^{\infty} \\
 & \ = \ R \frac{R^{-1+(2\log R)^{-\vareps}}}{1-(2\log R)^{-\vareps}} 
 \ = \ \frac{R^{(2\log R)^{-\vareps}}}{1-(2\log R)^{-\vareps}}.
\end{align*}
Similar to the previous calculus, 
we find that $\log R^{(2\log R)^{-\vareps}} \ = \  2^{-\vareps} (\log R)^{1-\vareps}$. Then 
it follows that 
$$
\limsup_{R\to\infty} R \int_R^{\infty} u^{-1-\alpha(u)}du<\infty.
$$
Thus the process is recurrent for $\vareps\ge 1$.

Now we show that the process is transient if $0<\vareps<1$.  In order to show this, 
recall that the characteristic function $\varphi$ of the process is defined by
$$
\varphi(\xi)  \ =  \ \int_{\mathbb R} \bigl(1-\cos (\xi h) \bigr) 
|h|^{-1-\alpha(|h|)}dh, \quad \xi \in {\mathbb R}
$$
and the process is recurrent if and only if for some (or equivalently,  for all) $r>0$, 
$$
\int_{\{ |\xi|<r\}} \frac {d\xi}{\varphi(\xi)} =\infty.
$$
(see \cite{S99}). Then we will prove that 
$\int_{\{ |\xi|<r\}} \frac {d\xi}{\varphi(\xi)}<\infty$ for some $0<r\le 1$ provided that $0<\vareps<1$.
This means it is enough for us to estimate the function $\varphi$ on 
$\{\xi \in{\mathbb R}: \  |\xi|<1\}$.

Since $\varphi(0)=0$, we only consider the case $0<|\xi|<1$.

First assume that $0<\xi<1$. Then 
\begin{align*}
\varphi(\xi) & \ =  \ \int_{\mathbb R} \bigl(1-\cos (\xi h) \bigr) 
|h|^{-1-\alpha(|h|)}dh  \\
& \ \ge \ \int_{\{ \pi/2 <\xi x <\pi \}} \bigl(1-\cos (\xi h) \bigr) 
|h|^{-2+ (\log (|h|+e^2))^{-\vareps}}dh \\
& \ = \ \int_{\{ \pi/2 <u <\pi \}} \bigl(1-\cos u \bigr) 
\Bigl|\frac u{\xi}\Bigr|^{-2+ (\log (|u/\xi|+e^2))^{-\vareps}}\frac{du}{\xi} 
\quad \quad (\xi h=u) \\
& \ \ge \ \xi^{1-(\log(\pi/\xi +e^2))^{-\vareps}} 
\int_{\{ \pi/2 <u <\pi \}} \bigl(1-\cos u \bigr)  u^{-2}du \ \ge \ 
c \, \xi^{1-(\log(\pi/\xi +e^2))^{-\vareps}},
\end{align*}
where $c$ is a constant independent of $\xi$. Similarly, we can get a similar  
bound for $-1<\xi<0$: 
\begin{align*}
\varphi(\xi) & \ =  \ \int_{\mathbb R} \bigl(1-\cos (\xi h) \bigr) 
|h|^{-1-\alpha(|h|)}dh  \\
& \ \ge \ \int_{\{ -\pi/2 >\xi x >-\pi \}} \bigl(1-\cos (\xi h) \bigr) 
|h|^{-1-\alpha(|h|)}dh \\
& \ \ge \ (-\xi)^{1-(\log(\pi/(-\xi) +e^2))^{-\vareps}} 
\int_{\{ \pi/2 <u <\pi \}} \bigl(1-\cos u \bigr)  u^{-2}du \ \ge \ 
c\, (-\xi)^{1-(\log(\pi/(-\xi) +e^2))^{-\vareps}}.
\end{align*}
Thus it follows that 
$$
\varphi(\xi) \ge c |\xi|^{1-(\log(\pi/|\xi| +e^2))^{-\vareps}}, \quad 0<|\xi|<1.
$$
Then, noting $0<\vareps<1$, we find that 
\begin{align*}
\int_{B(1)} \frac{ d\xi}{\varphi(\xi)} 
 & \ \le \ c \int_{|\xi|<1}  |\xi|^{-1+(\log(\pi/|\xi| +e^2))^{-\vareps}} d\xi  \\
 & \ \le \ 2c \int_0^1 u^{-1+(\log (\pi/u+e^2))^{-\vareps}} du 
 \quad \quad  (\pi/u+e^2=t) \\
 &  \ = \ \pi c \int_{\pi+e^2}^{\infty} 
\bigl( \pi(t-e^2)^{-1} \bigr)^{-1+(\log t)^{-\vareps}} \cdot (t-e^2)^{-2}dt 
 \\
 & \ \le \ c' \int_{\pi+e^2}^{\infty} (t-e^2)^{-1-(\log t)^{-\vareps}}dt 
\quad \quad (\log t=s) \\
 & \ \le \ c' \int_{\log(\pi+e^2)}^{\infty} (e^s-e^2)^{-1-s^{-\vareps}} 
\cdot e^sdt  \ \le \ c'' \int_2^{\infty} e^{-s^{1-\vareps}}ds 
\quad \quad (s^{1-\vareps}=x)  \\
 & \ = \ \frac{c''}{1-\vareps} \int_{2^{1-\vareps}}^{\infty} e^{-x} 
x^{\vareps/(1-\vareps)}dx  \ \le \ c''' \, \Gamma\Bigl(\frac{\vareps}{1-\vareps}+1\Bigr) <\infty.
\end{align*}
Therefore the form/process is transient for $0<\vareps<1$.
\hfill \fbox{}


\begin{flushright}
{\hfill Kohei Suzuki
\\ \hfill Department of Mathematics, Faculty of Science
\\ \hfill Kyoto University
    \\ \hfill Sakyo-Ku, Kyoto, 606-8502, Japan}
\bigskip \\
{\hfill Toshihiro Uemura
\\ \hfill Department of Mathematics, Faculty of Engineering Science
\\ \hfill  Kansai University
    \\ \hfill Suita, Osaka, 564-8680 Japan}
  \end{flushright}

\end{document}